# NOTE ON THE ZEROS OF A DIRICHLET FUNCTION

**Ferry**, L.[(1)], **Ghisa**, D.[(2)] **and Muscutar**, F. A.[(3)]

**Abstract**: The existence of non trivial zeros off the critical line for a function obtained by analytic continuation of a particular Dirichlet series is studied. Our findings are in contradiction with some results of computations which were present in the field for a long time. In the first part of the note we illustrate how the approximation errors may have had as effect inexact conclusions and in the second part we prove rigorously our point of view.

**AMS** subject classification: 30C35, 11M26

Examples of functions obtained by analytic continuation of Dirichlet series, which have zeros off the line $\operatorname{Re} s = 1/2$ are important for the purpose of circumscribing the field where the Riemann Hypothesis might be true. Such a candidate is attributed by some mathematicians to Davenport and Heilbronn (1936) (see [1]) and by others to Titchmarsh [5] (see [4]). In [6], R. C. Vaughan provided an elementary clear presentation of the respective example using the robust theory from [3]. He starts with that Dirichlet character $\chi$ modulo 5 for which $\chi(2) = i$. For this character

(1)  $\tau(\chi) = \sum_{k=1}^{5} \chi(k) e(\frac{k}{5}) = 2(-\sin\frac{\pi}{5} + i\sin\frac{2\pi}{5})$ and

(2)  $\epsilon(\chi) = \frac{\tau(\chi)}{i\sqrt{5}} = \frac{1}{\sqrt{2}}(\sqrt{1 + \frac{1}{\sqrt{5}}} + i\sqrt{1 - \frac{1}{\sqrt{5}}})$ and $|\epsilon(\chi)| = 1$.

Repeating the same computation for the character $\overline{\chi}$ it can be easily found that $\epsilon(\overline{\chi}) = \overline{\epsilon(\chi)}$. With the notation used in *Mathematica* implementation of Dirichlet L-functions, $L(5,2,s) = \sum_{n=1}^{\infty} \chi(n) n^{-s}$ and $L(5,3,s) = \sum_{n=1}^{\infty} \overline{\chi}(n) n^{-s}$ for $\operatorname{Re} s > 1$, it is known that $L(5,2,s)$ and $L(5,3,s)$ have analytic continuations to the whole complex plane and they are entire functions. Moreover, by corollary 10.9 of [3], they verify the functional equations:

(3)  $L(5,2,s) = \epsilon(\chi) 2^s \pi^{s-1} 5^{1/2-s} \Gamma(1-s)(\cos\frac{\pi s}{2}) L(5,3,1-s)$

(4)  $L(5,3,s) = \overline{\epsilon(\chi)} 2^s \pi^{s-1} 5^{1/2-s} \Gamma(1-s)(\cos\frac{\pi s}{2}) L(5,2,1-s)$

For $e^{2i\theta} = \epsilon(\chi)$ an elementary computation gives

(5)  $\tan^2\theta = [\sqrt{2} - \sqrt{1 + \frac{1}{\sqrt{5}}}] / [\sqrt{2} + \sqrt{1 + \frac{1}{\sqrt{5}}}]$, hence $\tan\theta = 0.284079\ldots$

For this value of $\tan\theta$ and taking into account [6], the function

(6)  $f(s) = \frac{1}{2}\sec\theta[e^{-i\theta}L(5,2,s) + e^{i\theta}L(5,3,s)] =$

$\frac{1}{2}\{[L(5,3,s) + L(5,2,s)] + i\tan\theta[L(5,3,s) - L(5,2,s)]\} =$

$\frac{1}{1^s} + \frac{\tan\theta}{2^s} - \frac{\tan\theta}{3^s} - \frac{1}{4^s} + \frac{1}{6^s} + \ldots$ verifies the functional equation:

(7)  $f(s) = 2^s\pi^{s-1}5^{1/2-s}\Gamma(1-s)(\cos\frac{\pi s}{2})f(1-s)$

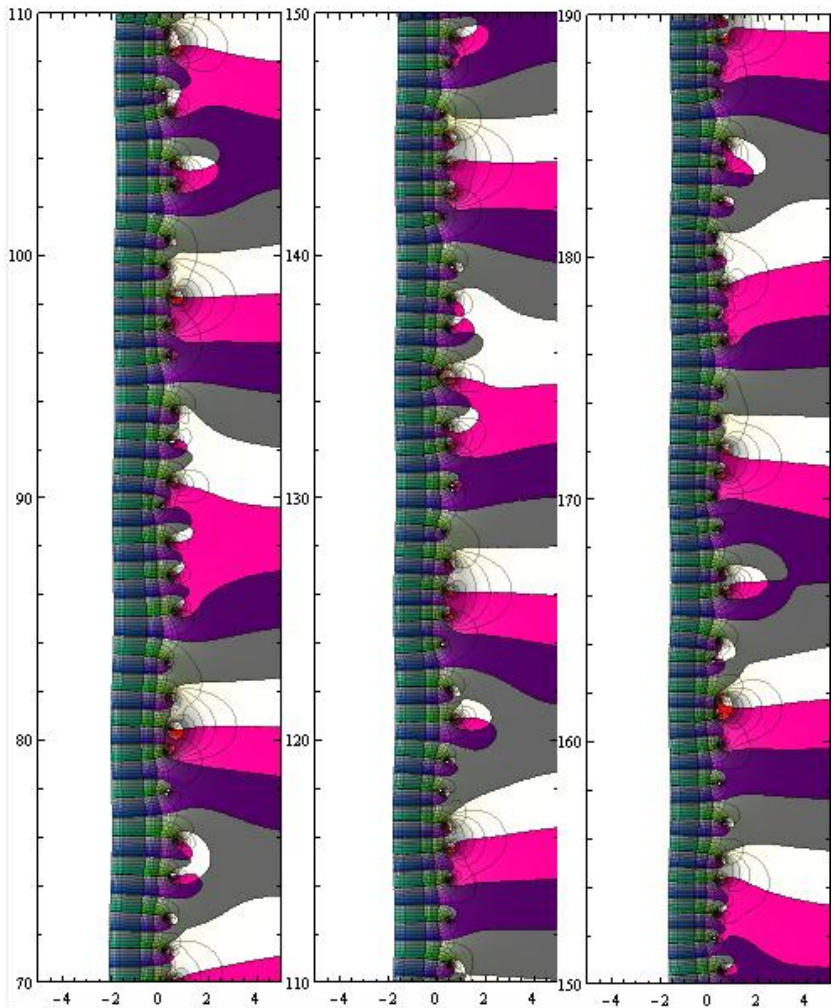

Fig. 1

We notice that the multiplier of $f(1 - s)$ can cancel only for real values of $s$. Therefore by (7), for $t \neq 0$, we have that $f(\sigma + it) = 0$ if and only if $f(1 - \sigma + it) = 0$, given the fact that, by (6), $f(1 - s) = \overline{f(1 - \bar{s})}$. In analogy with the case of the Riemann Zeta function, we call the non real zeros of $f(s)$ non trivial zeros. The previous remark tells us that the non trivial zeros of $f(s)$ appear in quadruplets two by two symmetric with respect to the real axis and two by two symmetric with respect to the critical line $\operatorname{Re} s = 1/2$.

Some computation done in [4] suggested that $.808517 + 85.699348i$, $0650830 + 114.163343i$, $.574356 + 166.479306i$ and $.724258 + 176.702461i$ are zeros of the function $f(s)$, which are obviously off the critical line. There was no indication that the symmetric points with respect to the critical line are also zeros of $f(s)$, as they should be, by (6) and (7). In [1] a confirmation of zeros indicated in [4] is implied and some more off line zeros are assumed by using linear combinations of Dirichlet L-functions.

Those linear combinations are supposed to satisfy functional equations similar to (7). No proof is given for this last affirmation. We have seen previously that this happened for functions defined by conjugate Dirichlet characters and only for one particular value of the parameter, fact which makes doubtful the existence of such an equation in general.

Fig. 1 above suggests that for the function $f(s)$ there are no off critical line non trivial zeros in the interval $t \in [80, 180]$ considered in [4] and [1], despite of the obvious non alignment of the zeros. This last phenomenon is certainly due to the errors of approximation of $\tan \theta$.

Indeed, since the series (6) converges extremely slowly in the neighborhood of a zero, the accumulation of the respective errors can have as effect a slight perturbation of the location of those zeros.

However, by (7) these zeros should be on the critical line, since otherwise they must be accompanied by symmetric zeros with respect to the critical line. In other words, a hint of their symmetry should have been persisted after perturbation, in the sense that some couple of zeros should have had close values for the $t$-coordinate. There is not evidence of such a symmetry in the relevant interval.

What one can see in the picture is a family of orthogonal curves which are mapped by $f(s)$ onto the orthogonal net formed with circles centered at the origin and rays issuing from the origin. The net is implemented in *Mathematica* and consists of 12 rays, two consecutive ones making an angle of $\pi/6$ between them and a number of circles centered at the origin, the radii of which are 1/12, 2/12,..., 12/12, 2, 10, $10^2$,…

The pre-image of circles helps locate the zeros in an obvious way and at the same time gives a clear description of the mapping at different orders of magnitude of the rays.

They also give an indication about the relative *closeness* of the zeros. The chosen interval for $t$ contains the four points indicated in [4] as zeros of $f(s)$. It also exhibits strips $S_k$ with $t$ in a span averaging 9 units (indeed in the interval [77,185] there are exactly 12 strips).

Every strip is mapped by $f(s)$ onto the whole complex plane with a slit alongside the interval $[1,+\infty]$ of the real axis. The mapping is $j_k$ to one where $j_k$ is the number of zeros of $f(s)$ in $S_k$. This number coincides with that of the components of the pre-image of small circles centered at the origin.

If the points in [4] were zeros of the function, then they should coincide or should be very close to some of these perturbed zeros. It can be noticed that the point .574356+166.479306i signaled in [4] is located somewhere between the two zeros belonging to a bounded component of the pre-image of the unit disc. The two zeros can be considered close to each other, since they belong to a component of the pre-image of the disc of radius 1/2.

This can be seen by noticing that there are six curves turning around both zeros. The other points indicated in [4] are too far from the zeros appearing in this picture in order to be confused with one of them. However, even this point cannot be taken as one of them since then the other should be the symmetric with it with respect to the critical line, which is not the case.

Thus, the points indicated in [4] as zeros are in disagreement with the functional equation (7).

Our analysis was prompted by the fact that in [2] it has been proved that, with a proper definition of the non trivial zeros, all the Dirichlet L-functions verify the Riemann Hypothesis. Moreover, a paper has been submitted to the attention of R. C. Vaughan in which the author affirms that the Riemann Hypothesis is a theorem for a large class of functions which verify some very general conditions and a functional equation similar to (7).

Doubts were raised about such a possibility having in view the results in [4] and [1].

We point out that those results are necessarily affected by errors of approximation, which are strictly localized. Our approach is global and allows one a view of the conformal mapping, and in particular to detect symmetric zeros if they existed. The method takes advantage of the exact implication of the equation (7) regarding the symmetry of those zeros.

An alternative approach has been developed by one of us[(*)] having at its core the comparison of partial sums by using the Java Applet *Critical Strip Explorer.* Applying it to a false zero provides a contradiction of the functional equation. Although this test is a rather experimental complement, it gave us a clue about non zeros of the points signaled in [4].

Finally, we notice that the purpose of the formula (6) was certainly to give an example of Dirichlet series whose analytic continuation does not verify the Riemann Hypothesis, although it verifies a functional equation, and this supposedly because the series is not an Euler product.

Opposite to this idea is our finding that the series does not need to be expressible as an Euler product in order for its continuation $f(s)$ to verify the Riemann Hypothesis. However, it needs to verify a functional equation relating the values of $f(s)$ and $\overline{f(1-\bar{s})}$ by the intermediate of a multiplier.

This is the case for the function (6) and in our view the off critical line points signaled in [4] and [1] are not true zeros. The true non trivial zeros of the function (6) must lie on the critical line, as shown in the following:

**Theorem**: *In every strip $S_k$ of $L(5,2,s)$ the number of non trivial zeros of $L(5,2.s)$ and $L(5,3,s)$ is the same. These zeros and those of $f(s)$ lie all on the critical line.*

*Proof:* Let $j_k$ be the number of non trivial zeros $s_{k,j} = 1/2 + it_{k,j}$ of $L(5,2,s)$ situated in $S_k$. This function maps every interval $[s_{k,j}, s_{k,j+1}]$ onto a closed curve passing through the origin. The same thing does the function (4), hence there is a point $s'_{k,j}$ such that $L(5,3,s'_{k,j}) = 0$. Choosing conveniently the allocation of zeros of $L(5,3,s)$ at the limit of two strips, we obtain the same numbers as for $L(5,2,s)$.

In a similar way we find that every strip $S'_j$ of $L(5,3,s)$ contains the same number of non trivial zeros of $L(5,3,s)$ and $L(5,2,s)$. The fundamental domains $\Delta_{k,j} \subset S_k$ of $L(5,2,s)$ are mapped conformally by $L(5,2,s)$ onto the complex plane with a slit. These domains are bounded by components of the pre-image by $L(5,2,s)$ of the interval $(1,+\infty)$ of the real axis and arcs connecting the points $u_{k,j}$ for which $L(5,2,u_{k,j}) = 1$ with zeros $v_{k,j}$ of $L'(5,2,s)$. Since $\lim_{\sigma \to +\infty} L(5,2,\sigma + it) = 1$, the point $s = \infty$ on the Riemann sphere is considered as one of $u_{k,j}$. The fundamental domains of $L(5,3,s)$ are also mapped conformally by $L(5,3,s)$ onto the complex plane with a slit. Thus, every zero $s_{k,j}$ of $L(5,2,s)$ belongs to two fundamental domains, $\Delta_{k,j}$ of $L(5,2,s)$ and say $\Delta'_{k,j}$ of $L(5,3,s)$. In Fig. 2 below, the corresponding curves $\Gamma_{k,j}$ and $\Upsilon_{k,j}$ appear as intertwining curves if $j \neq 0$. This might not be true, as visible in the respective figure if $j = 0$ or if the curves are $\Gamma'_k$ and $\Upsilon'_k$. However, every strip $S_k$ contains a unique curve $\Gamma_{k,0}$ and a unique curve $\Upsilon_{k,0}$ which implies that the number of zeros of the two functions situated in $S_k$ is the same.

By [2], the non trivial zeros $s_{k,j}$ of $L(5,2,s)$ and the non trivial zeros $s'_{k,j}$ of $L(5,3,s)$ lie all on the critical line, i.e. $s_{k,j} = 1/2 + it_{k,j}$ and $s'_{k,j} = 1/2 + it'_{k,j}$. If $s_{k,j} = s'_{k,j}$, then by (6), this is a zero of $f(s)$ and it belongs to the critical line. If not, let us notice that the image by $e^{-i\theta}L(5,2,s)$ of the interval $\{s \mid s = 1/2 + it,\ t_{k,j} \leq t \leq t_{k,j+1}\}$ is a closed curve passing through the origin.

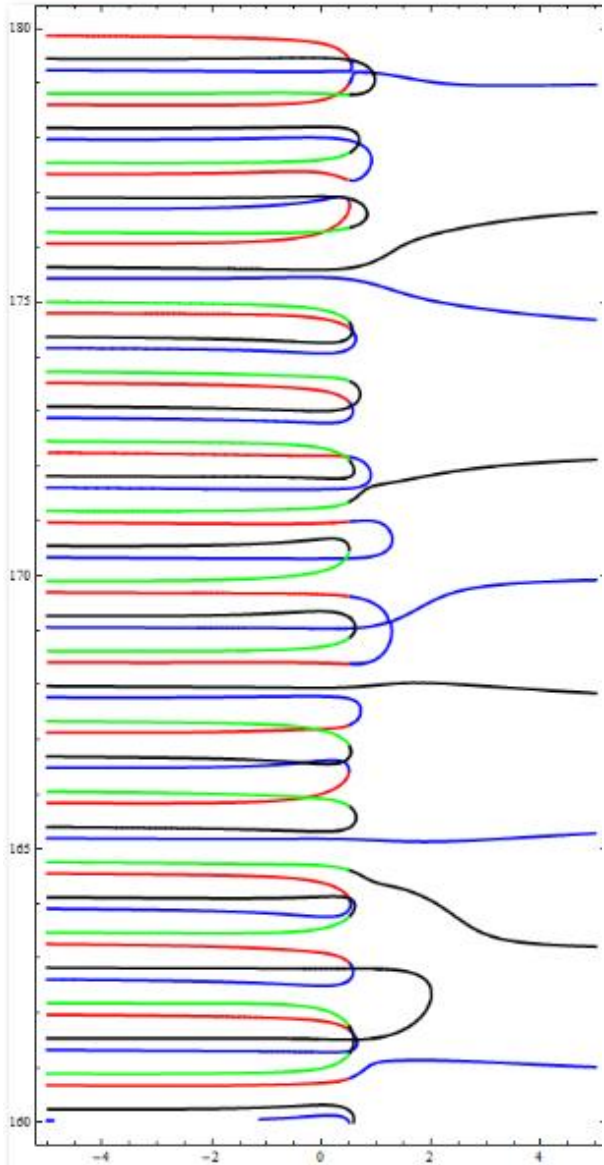

Fig. 2

The same is true for the image by $-e^{i\theta}L(5,3,s)$ of the interval $\{s \mid s = 1/2 + it, t'_{k,j} \le t \le t'_{k,j+1}\}$. The two curves may have in common only the origin, or they can intersect each other at another point corresponding to a value $t_0$ of $t$. In the last case $s_0 = 1/2 + it_0$ is by (6) a zero of $f(s)$, since $e^{-i\theta}L(5,2,s_0) = -e^{i\theta}L(5,3,s_0)$.

On the other hand, if for an $s_0 = \sigma_0 + it_0$ we have $f(s_0) = 0$, then by (7) this happens if and only if $f(1 - \sigma_0 + it_0) = 0$ and by (6) if and only if

(8) $\quad (1 - i\tan\theta)L(5,2,\sigma_0 + it_0) + (1 + i\tan\theta)L(5,3,\sigma_0 + it_0) = 0$, respectively

(9) $\quad (1 - i\tan\theta)L(5,2,1-\sigma_0+it_0) + (1+i\tan\theta)L(5,3,1-\sigma_0+it_0) = 0$

From (8) and (9) it can be easily found

(10) $\quad L(5,2,\sigma_0+it_0) / L(5,3,\sigma_0+it_0) = L(5,2,1-\sigma_0+it_0) / L(5,3,1-\sigma_0+it_0)$

which by (3), (4) and the main result in [2] is possible only if $\sigma_0 = 1/2$. Indeed, if $\sigma_0 = 1/2$ then (10) is true for any $t$, in particular for $t_{k,j} \leq t \leq t_{k,j+1}$, where $1/2.+it_{k,j}$ and $1/2 + it_{k,j+1}$ are consecutive zeros of $L(5,2,s)$. For $\sigma_0 \neq 1/2$ the two terms in (10) are on different sides of the image by $L(5,2,s) / L(5,3,s)$ of this segment, hence they are different. Finally $f(\sigma_0+it_0) = 0$ only if $\sigma_0 = 1/2$, which proves completely the theorem. Fig.2 above illustrates this situation for the strip $S_k$ corresponding to $t$ approximately between 165 and 175. This strip contains 8 zeros of each one of the two Dirichlet $L$-functions. None of the zeros of $f(s)$ coincides with $.574356 + 166.702461i$ signaled in [4] and [1]. Also, a part of the strip $S_{k+1}$ is visible in Fig. 2 which is supposed to contain the zero $.724258 + 176.702461i$. The true zero should be instead on the critical line between the first two curves from this strip.

**Acknowledgements**: The authors are thankful to R. C. Vaughan for inspiring us to produce this note.

[1]lesandad.ferry@btinternet.com  [2]dghisa@yorku.ca  [3]fmuscuta@lorainccc.edu

[*] Les Ferry,